\documentclass{amsart}
\usepackage{amsthm}
\usepackage{amsmath}
\usepackage{amssymb}
\usepackage[margin=0.75in]{geometry}
\usepackage{enumerate}
\usepackage{enumitem}
\usepackage{mathrsfs}
\usepackage{float}
\usepackage{graphicx}
\usepackage{mathtools}
\usepackage{hyperref}
\usepackage{caption}
\usepackage{subcaption}

\newtheorem{theorem}{Theorem}
\newtheorem{corollary}[theorem]{Corollary}

\newtheorem{lemma}[theorem]{Lemma}
\newtheorem{definition}[theorem]{Definition}
\newtheorem{example}[theorem]{Example}
\newtheorem{remark}[theorem]{Remark}
\numberwithin{theorem}{section}
\numberwithin{equation}{section}

\allowdisplaybreaks

\begin{document}

	\title[Rigidity of Angles]{Rigidity of angles in inner product spaces}


\author{Aniruddha Deshmukh}
\address{Department of Mathematics, Indian Institute of Technology Indore}
\curraddr{}
\email{aniruddha480@gmail.com}
\thanks{The first author wishes to thank the Prime Minister Research Fellowship (PMRF) number $2101706$, that has funded the research work.}

\author{Ashisha Kumar}
\address{Department of Mathematics, Indian Institute of Technology Indore}
\curraddr{}
\email{akumar@iiti.ac.in}
\thanks{}


\date{}

\dedicatory{}

\commby{}

\begin{abstract}
In this article, we discuss the equality of two inner products on a vector space. Particularly, we look at some geometric properties that are given to a vector space by an inner product namely, length and angle, and we ask under what conditions on these geometric properties do two inner products coincide. 
Some parts of this article can be found as exercise problems in various books on geometry (see for instance, \cite{LeeRM}). We give a slight generalization to these exercises and exposit the rigid nature of angles in an inner product space. 
\end{abstract}

\maketitle

	\section{Introduction}
	\label{IntroductionSection}
	Inner products form an important part of study of Linear Algebra, Functional Analysis and other fields of mathematics. Loosely speaking, they define ``\textit{angles}" and ``\textit{projections}" between vectors of a vector space. This gives tools of geometry to what was purely algebra (vector spaces) and also forms a connection to topology (normed spaces). Clearly, if two inner products are equal, the geometry defined by them will be the same, i.e., the angles between vectors will be the same with respect to the two inner products, the lengths of the vectors will be the same and hence the notion of distance with respect to the two inner products will be the same. However, when can we say that two inner products are equal? Are there any equivalent conditions that we can find in terms of geometric properties? This article deals with these questions, and answers them positively. The question of reconstructing an inner product solely from the knowledge of all the lengths is a common exercise problem, famously known as the \textit{polarization identities} (see for instance, \cite{Krishnamurthy}, \cite{HoffmannKunze}, or \cite{KumaresanLA}). One can also find the question of knowing an inner product from the knowledge of angles between \textbf{all} pairs of vectors (see for instance, \cite{LeeRM}, where the exercise is asked in the technicality of Riemannian geometry). This, in fact, has been our motivation to write this article.
	

	We give a natural generalization to these exercise problems. Particularly, the article is dedicated to determining the conformality of two inner products with the knowledge of only \textbf{one} fixed angle. We will make this statement precise in a moment. At this point, though, we would like to remark that although the question (and its answer) is not ground-breaking, we could not find any references speaking of it. This article also serves the purpose of archiving a small, yet important, observation about the rigidity of geometry\footnote{By the word ``\textit{geometry}", we mean to tools of measuring angles, lengths and distances. These tools are obtained with the help of inner products.} in inner product spaces.

The article is organized as follows: Section \ref{PreliminariesSection} gives basic definitions and results that we will be using throughout the article. We will also look at the question we want to pose in more detail and try to understand what geometric properties could give us the equivalence we desire. Section \ref{RealCase} deals with the real vector spaces, where we prove our main results and answer the question raised at the beginning of this article positively. Section \ref{ComplexCase} then extends these results to complex vector spaces. Most of the results of this section are analogous to that of real vector space setting, and we therefore omit the detailed proofs. 
Section \ref{ConclusionSection}, which is the last section of this article, includes some concluding remarks, and an exploratory comment about conformality of Riemannian metrics. We pose a question concerning the same at the end of this article.
	\section{Preliminaries}
	\label{PreliminariesSection}
	Throughout this article, we assume that the reader is thorough with the basic definitions of a field and a vector space (a few references are \cite{HoffmannKunze}, \cite{Krishnamurthy}, and \cite{KumaresanLA}). 
	The simplest examples of vector spaces are $\mathbb{R}^n$ and $\mathbb{C}^n$, the $n$-tuples of real and complex numbers with component wise addition and scalar multiplication, over $\mathbb{R}$ and $\mathbb{C}$, respectively. There are many more examples of vector spaces, of which sequence spaces and function spaces are often studied and find their use in a variety of fields of mathematics (see \cite{RudinRA} and \cite{RudinRCA}, where these spaces are discussed from many viewpoints). In this article, we do not go into the details of different vector spaces, but only look at the abstract spaces over $\mathbb{R}$ and $\mathbb{C}$. 

	The study of mathematics at an early stage involves a lot of plane (Euclidean) geometry. This mostly involves angles and lengths of line segments. We can then try to define angles (and lengths) between vectors of an arbitrary vector space.
	Indeed, the inherent geometric idea comes from the case of $\mathbb{R}^2$, where the angle between $x, y \in \mathbb{R}^2$ is given by
	\begin{equation}
		\label{AngleR2}
		\cos \theta = \frac{x_1y_1 + x_2y_2}{\sqrt{x_1^2 + x_2^2} \sqrt{y_1^2 + y_2^2}}.
	\end{equation}
	For a detailed treatment of the topic, we refer the reader to \cite{KumaresanLA} and the beautifully written lecture notes at: \\
	{\small{\center{\url{www.mit.edu/~hlb/StantonGrant/18.02/details/tex/lec1snip2-dotprod.pdf.}}}}\\
	
	The motivation to define inner products comes from this observation in the Euclidean geometry. We now give the formal definition of inner products, and hence angles.
	\begin{definition}[Inner Product]
		Let $V$ be a vector space over a field $\mathbb{F}$\footnote{Here, our field $\mathbb{F}$ is either $\mathbb{R}$ or $\mathbb{C}$.}. An inner product is a function $\langle \cdot, \cdot \rangle : V \times V \rightarrow \mathbb{F}$ that satisfies:
		\begin{enumerate}
			\item[(IP1)] For every $x \in V$, we have $\langle x, x \rangle \geq 0$ and $\langle x, x \rangle = 0$ if and only if $x = 0$. \hfill [Positive Definiteness]
			\item[(IP2)] For every $x, y \in V$, we have $\langle x, y \rangle = \overline{\langle y, x \rangle}$. \hfill [Conjugate symmetry]
			\item[(IP3)] For every $x, y, z \in V$ and $\alpha, \beta \in \mathbb{F}$, we have $\langle \alpha x + \beta y, z \rangle = \alpha \langle x, z \rangle + \beta \langle y, z \rangle$. \hfill [Sesquilinearity\footnote{Notice that (IP2) and (IP3) together give us $\langle x, \alpha y + \beta z \rangle = \overline{\alpha} \langle x, y \rangle + \overline{\beta} \langle x, z \rangle$. That is, in the second input of inner product, addition behaves well but the scaling doesn't. This is the reason (IP3) is called ``\textit{sesquilinearity}", which means that the inner product as a function of two variables is ``\textit{one-and-a-half linear}". Had it also respected scaling in the second input, we would have called it bilinear.}]
		\end{enumerate}
	\end{definition}
	\begin{remark}
		\normalfont
		If $\mathbb{F} = \mathbb{R}$, then the conjugate symmetry becomes symmetry, i.e., $\langle x, y \rangle = \langle y, x \rangle$, and the sesquilinearity becomes bilinearity. The case of real inner product spaces is, therefore, easier to handle.
	\end{remark}
	One of the natural examples of inner products 
		is the ``\textit{usual}" or ``\textit{standard}" or ``\textit{Euclidean}" inner product on $\mathbb{R}^n$ and $\mathbb{C}^n$. It is given by 
		$\langle z, w \rangle = \sum\limits_{i = 1}^{n} z_i \overline{w_i},$
		where $z = \left( z_1, z_2, \cdots, z_n \right), w = \left( w_1, w_2, \cdots, w_n \right) \in \mathbb{C}^n$. \\
		With this understanding, we see that for any $x \in V$, $\| x \| = \sqrt{\langle x, x \rangle}$ gives a norm on $V$. The proof of this fact uses the famous Cauchy-Schwartz inequality whose proof can be found in \cite{HoffmannKunze}, \cite{Krishnamurthy}, or any other standard linear algebra text-book.
		\begin{theorem}[Cauchy-Schwartz]
			\label{CSI}
			Let $\left( V, \langle \cdot, \cdot \rangle \right)$ be an inner product space. Then, for any $x, y \in V$, we have
			\begin{equation}
				\label{CauchySchwartzInequality}
				\left| \langle x, y \rangle \right| \leq \| x \| \| y \|.
			\end{equation}
		\end{theorem}
		The Cauchy-Schwartz inequality \eqref{CauchySchwartzInequality} is remarkable in the sense that with this, we can now say (at least for real vector spaces) that
		$$-1 \leq \dfrac{\langle x, y \rangle}{\| x \| \| y \|} \leq 1,$$
		when $x \neq 0$ and $y \neq 0$. Thus, there is a \textit{unique} real number $\theta \in \left[ 0, \pi \right]$ such that $\cos \theta = \frac{\langle x, y \rangle}{\| x \| \| y \|}$. We define this number $\theta$ as the \textit{\textbf{angle}} between $x$ and $y$.

		In the setting of complex inner product spaces, there are various ways to define the notion of angles. In fact, one can observe that a lot of linear algebra textbooks avoid this subject, probably to reduce confusion among its readers. The problem is that if we start with a complex inner product space $\left( V, \langle \cdot, \cdot \rangle \right)$, then for any (non-zero) vectors $x, y \in V$, the quantity $\langle x, y \rangle$ is, in general, a complex number. Hence, the quantity $\frac{\langle x, y \rangle}{\| x \| \| y \|}$, in general, cannot be the cosine of any \textit{real} number! Nonetheless, we can define it to be the cosine of a complex number, say $\theta_{C}$, and call it the \textit{\textbf{complex angle}} between $x$ and $y$. However this definition has no geometric intuition behind it. We have simply borrowed the formula we know in a real vector space!

		One can then ask if we can use $\left| \langle x, y \rangle \right|$, which is a real quantity and the Cauchy-Schwartz inequality guarantees that $\frac{\left| \langle x, y \rangle \right|}{\| x \| \| y \|}$ is the cosine of some real number. However, the quantity is always non-negative! Therefore, if we define the angle using the aforementioned quantity, we cannot talk about ``obtuse" angles. This seems counter-intuitive since geometrically, we would at least want $x$ and $-x$ to face in ``opposite directions", and hence have an angle $\pi$ between them. 

		To understand the ``natural" choice\footnote{By the word ``natural", we mean a more geometrically intuitive choice of defining angles.} of expression to define angles in a complex vector space setting, we will have to work a bit. While in $\mathbb{R}^2$, the (standard) inner product had a geometric meaning and this meaning was closely related to the angle between given vectors, this is neither intuitive nor straight-forward in complex vector spaces. Let us start with the simplest complex vector space $\mathbb{C}$.  

\begin{figure}[ht!]
			\centering
			\includegraphics[scale=4.0]{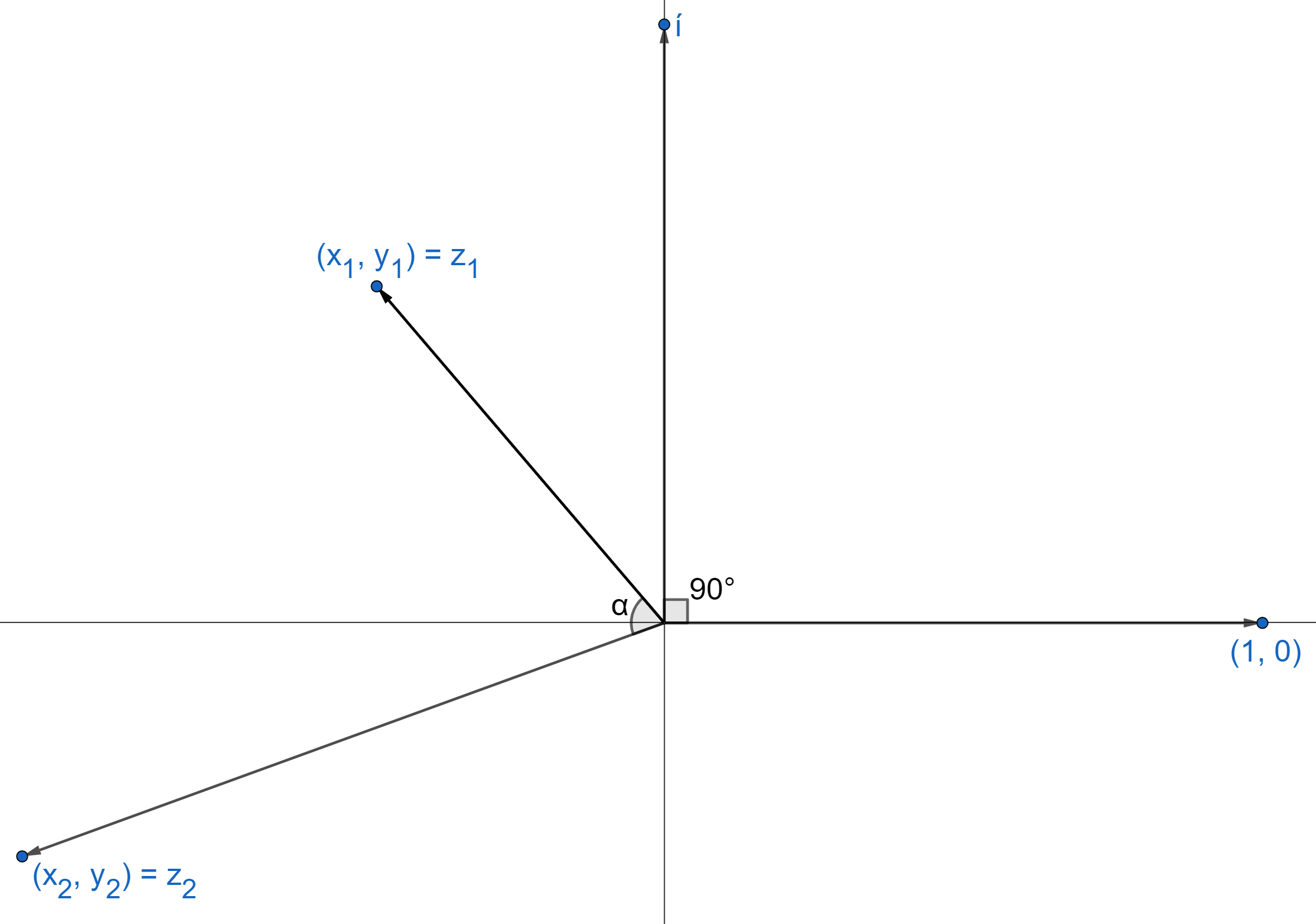}
			\caption{The identification of $\mathbb{C}$ with $\mathbb{R}^2$.}
			\label{Identification}
		\end{figure}

		Let us first observe that $\mathbb{C}$ is also a real vector space and with the real structure, it is isomorphic to $\mathbb{R}^2$. A complex number $z = a + \iota b$ can be identified with the tuple $\left( a, b \right) \in \mathbb{R}^2$. The identification is shown in Figure \ref{Identification}. 
		Both the spaces have inner products defined on them. On $\mathbb{C}$, we have $\langle z_1, z_2 \rangle_{\mathbb{C}} = z_1 \overline{z_2}$, while on $\mathbb{R}^2$, we have $\langle \left( x_1, x_2 \right), \left( y_1, y_2 \right) \rangle_{\mathbb{R}^2} = x_1 y_1 + x_2 y_2$. 
		Let us start with $z_1, z_2 \in \mathbb{C}$, and let the corresponding vectors in $\mathbb{R}^2$ be $\left( a_1, b_1 \right), \left( a_2, b_2 \right)$. We first see that
		$$\langle z_1, z_2 \rangle_{\mathbb{C}} = z_1 \overline{z_2} = \left( a_1 a_2 + b_1 b_2 \right) - \iota \left( a_1 b_2 - a_2 b_1 \right).$$
		That is, $\text{Re} \left( \langle z_1, z_2 \rangle_{\mathbb{C}} \right) = a_1 a_2 + b_1 b_2 = \langle \left( a_1, b_1 \right), \left( a_2, b_2 \right) \rangle_{\mathbb{R}^2}$.
		Coming to $\mathbb{R}^2$, we know that the angle between the vectors $\left( a_1, b_1 \right)$ and $\left( a_2, b_2 \right)$ is given by $\cos \theta = \dfrac{\langle \left( a_1, b_1 \right), \left( a_2, b_2 \right) \rangle_{\mathbb{R}^2}}{\| \left( a_1, b_1 \right) \|_{\mathbb{R}^2} \| \left( a_2, b_2 \right) \|_{\mathbb{R}^2}}$. However, with the identification that we have made, it is easy to see that $\| z \|_{\mathbb{C}} = \| \left( a, b \right) \|_{\mathbb{R}^2}$. Thus, we have
		$$\cos \theta = \dfrac{\text{Re} \left( \langle z_1, z_2 \rangle_{\mathbb{C}} \right)}{\| z_1 \|_{\mathbb{C}} \| z_2 \|_{\mathbb{C}}}.$$
		We define this $\theta$ as the angle between the vectors $z_1$ and $z_2$ in the complex inner product space $\mathbb{C}$. Abstracting this observation, we define the angle between two vectors $x$ and $y$ in a complex inner product space $V$ as the unique real number $\theta \in \left[ 0, \pi \right]$ such that
		$$\cos \theta = \dfrac{\text{Re } \langle x, y \rangle}{\| x \| \| y \|}.$$
		\begin{remark}
			\normalfont
			Since we have
			$\left| \text{Re} \left( \langle x, y \rangle \right) \right| \leq \left| \langle x, y \rangle \right| \leq \| x \| \| y \|,$ 
			$\theta \in \left[ 0, \pi \right]$ is well-defined.
		\end{remark}
		\begin{remark}
			\label{DifferentDefinitionsComplexAngle}
			\normalfont
			We would like to list various ways to define angles between vectors of a complex inner product space $\left( V, \langle \cdot, \cdot \rangle \right)$.
			\begin{equation}
				\label{GeometricAngle}
				\cos \theta = \frac{\text{Re } \langle x, y \rangle}{\| x \| \| y \|},
			\end{equation}
			\begin{equation}
				\label{ComplexAngle}
				\cos \theta_C = \frac{\langle x, y \rangle}{\| x \| \| y \|},
			\end{equation}
			\begin{equation}
				\label{HermitianAngle}
				\cos \theta_H = \frac{\left| \langle x, y \rangle \right|}{\| x \| \| y \|},
			\end{equation}
			\begin{equation}
				\label{PsuedoAngle}
				\varphi = \text{Arg } \cos \theta_C,
			\end{equation}
			\begin{equation}
				\label{KahlerAngle}
				\cos \theta_K \sin \theta = \frac{\text{Re } \langle \iota x, y \rangle}{\| x \| \| y \|}.
			\end{equation}
			All these definitions arise in various aspects of geometry and (complex) analysis. Indeed, Equation \eqref{GeometricAngle} is the geometrically intuitive definition of angles. It is called the \textbf{\textit{Euclidean angle}} between vectors. A few Linear Algebra textbooks (see for instance, \cite{Golan}) use this definition. Equation \eqref{HermitianAngle} is called the \textbf{\textit{Hermitian angle}} and can be found in \cite{Ricardo}, \cite{Coolidge}, and many other journal references (see the references of \cite{Scharnhorst}). The geometric meaning of the Hermitian angle is given in \cite{Scharnhorst}. However, it does not go well with the intuition that $\mathbb{C}$ can be (canonically) identified with $\mathbb{R}^2$. The notion of (Kasner's) \textbf{\textit{pseudo-angle}} given by Equation \eqref{PsuedoAngle} is used by Kasner in \cite{Kasner1} in the context of conformality of functions of (several) complex variables. The \textbf{\textit{K\"{a}hler angle}} given by Equation \eqref{KahlerAngle} is studied when we look at the multiplication by $\iota$ as a linear map on the real vector space\footnote{By this, we mean that if $V$ is a given complex vector space, we can also look at it as a real vector space, where $\left\lbrace x, \iota x \right\rbrace$ is linearly independent.} $V$. The K\"{a}hler angle measures how far is a (real) two-dimensional plane from being invariant under the action of $\iota$. For a treatment on K\"{a}hler angles, we refer the reader to \cite{Scharnhorst} and the references therein.

			For all considerations of this article, we will only work with the Euclidean (geometric) angle given by Equation \eqref{GeometricAngle}.
		\end{remark}
		\begin{remark}
			\normalfont
			Notice that the angle between $x$ and $y$ will not be defined when $x = 0$ or $y = 0$. This is because in this case we have $\langle x, y \rangle = 0$\footnote{This can be easily verified using the sesquilinearity of the inner product and the properties of the zero vector.} and one of $\| x \|$ or $\| y \|$ is zero. Hence, a \textbf{unique} $\theta$ cannot be determined in this case for which $\text{Re } \langle x, y \rangle = \| x \| \| y \| \cos \theta$.
		\end{remark}
		An important aspect in the main result of this article would be ``perpendicular" vectors. Closely related to being perpendicular is the notion of being orthogonal, which we now briefly view.
		\begin{definition}[Orthogonality]
			\label{Orthogonality}
			Let $\left( V, \langle \cdot, \cdot \rangle \right)$ be an inner product space. Two vectors $x, y \in V$ are orthogonal if $\langle x, y \rangle = 0$.
		\end{definition}
		Clearly, in the case of real inner product spaces, we see that two vectors $x, y \in V$ are orthogonal if and only if the angle between them is $\frac{\pi}{2}$. However, since the definition of angle between vectors in complex inner product spaces involves only the real part of the inner product, such an equivalence cannot be expected. Surely, orthogonal vectors have $\frac{\pi}{2}$ angle between them but the converse is not true.
		\begin{example}
			\label{PerpendicularNotOrthogonal}
			\normalfont
			Let us consider the complex vector space $\mathbb{C}$ with the standard inner product given by $\langle z, w \rangle = z \overline{w}$, for $z, w \in \mathbb{C}$. With this inner product, we see that $\langle 1, \iota \rangle = - \iota$, so that $\text{Re } \langle 1, \iota \rangle = 0$. Thus, the angle between $1$ and $\iota$ is $\frac{\pi}{2}$ (see Figure \ref{Identification}), but clearly they are not orthogonal.
		\end{example}
		Now that we have made the notion of angles in arbitrary vector spaces precise, we can ask the following question: \\
		\\
			\textbf{Question}: Given two inner products on a vector space $V$, say $\langle \cdot, \cdot \rangle_1$ and $\langle \cdot, \cdot \rangle_2$, under what conditions can these be equal?\\
			\\
		An easy (and correct) answer to this is when for every pair of vectors $x, y \in V$, we have $\langle x, y \rangle_1 = \langle x, y \rangle_2$. This follows from the definition of the equality of two functions! However, checking the values of these inner products at \textbf{every} pair of vectors is too much to ask for!

		Let us try to reduce this work. Since inner products are related to lengths of vectors, let us ask if the two inner products are same given that the lengths corresponding to these inner products are same. That is, we want to see whether $\langle \cdot, \cdot \rangle_1 = \langle \cdot, \cdot \rangle_2$ when for each $x \in V$, we have $\| x \|_1 = \| x \|_2$. Here, $\| \cdot \|_1$ denotes the norm corresponding to $\langle \cdot, \cdot \rangle_1$, and likewise for $\| \cdot \|_2$. This is in fact true, and is an immediate consequence of the polarization identities (see, for instance, \cite{Krishnamurthy}, \cite{HoffmannKunze}, or \cite{KumaresanLA}, where they appear as exercise problems). The proof immediately follows by using the definition of norms induced by inner products (i.e., $\| x \|^2 = \langle x, x \rangle$) and the axiomatic properties (IP1)-(IP3) of inner products. 
		\begin{theorem}[Polarization identities]
			\label{PITheorem}
			Let $V$ be a complex inner product space with inner product $\langle \cdot, \cdot \rangle$. Then, for every $x, y \in V$, we have
			\begin{equation}
				\label{ComplexPI}
				\langle x, y \rangle = \dfrac{1}{4} \left( \| x + y \|^2 - \| x - y \|^2 \right) + \dfrac{\iota}{4} \left( \| x + \iota y \|^2 - \| x - \iota y \|^2 \right).
			\end{equation}
			For a real inner product space $V$, we have for every $x, y \in V$,
			\begin{equation}
				\label{RealPI}
				\langle x, y \rangle = \dfrac{1}{4} \left( \| x + y \|^2 - \| x - y \|^2 \right).
			\end{equation}
		\end{theorem}
		From Theorem \ref{PITheorem}, it is clear that if $\| \cdot \|_1 = \| \cdot \|_2$, then $\langle \cdot, \cdot \rangle_1 = \langle \cdot, \cdot \rangle_2$. An immediate corollary of the polarization identities will be beneficial to us for our results.
		\begin{corollary}
			\label{ConformalCorollary}
			Let $V$ be a (real or complex) vector space with inner products $\langle \cdot, \cdot \rangle_1$ and $\langle \cdot, \cdot \rangle_2$. Let $\| \cdot \|_1$ and $\| \cdot \|_2$ be the corresponding norms. If there is some $c > 0$ such that $\| \cdot \|_1 = c \| \cdot \|_2$, then we have $\langle \cdot, \cdot \rangle_1 = c^2 \langle \cdot, \cdot \rangle_2$.
		\end{corollary}
		We have seen that angles are intimately connected with inner products, and in fact, are the main theme of this exposition. Therefore, we ask the question: If $\langle \cdot, \cdot \rangle_1$ and $\langle \cdot, \cdot \rangle_2$ define the same angles, can we conclude that $\langle \cdot, \cdot \rangle_1 = \langle \cdot, \cdot \rangle_2$? By ``define the same angles" we mean that for any $x, y \in V \setminus \left\lbrace 0 \right\rbrace$, the angle between $x$ and $y$ with respect to $\langle \cdot, \cdot \rangle_1$ is $\theta$ if and only if the angle between them with respect to $\langle \cdot, \cdot \rangle_2$ is $\theta$.\\
		\\
		The following example shows that this is indeed not the case.
	\begin{example}
		\label{SameAngleDifferentInnerProducts}
		\normalfont
		Let us consider $\langle \cdot, \cdot \rangle_1$, the standard inner product on $\mathbb{R}^2$, and $\langle \cdot, \cdot \rangle_2 = 2 \langle \cdot, \cdot \rangle_1$. For any $x, y \in \mathbb{R}^2 \setminus \left\lbrace 0 \right\rbrace$, let $\theta_1$ be the angle between $x$ and $y$ as defined by $\langle \cdot, \cdot \rangle_1$ and $\theta_2$ be the angle between $x$ and $y$ as defined by $\langle \cdot, \cdot \rangle_2$. Then,
		\begin{align*}
			\cos \theta_2 &= \dfrac{\langle x, y \rangle_2}{\| x \|_2 \| y \|_2} 
			= \dfrac{2 \langle x, y \rangle_1}{\sqrt{2 \langle x, x \rangle_1} \sqrt{2 \langle y, y \rangle_1}}
			= \dfrac{\langle x, y \rangle_1}{\| x \|_1 \| y \|_1}
			= \cos \theta_1.
		\end{align*}
		That is, the angles defined by the two inner products are the same. However, the two inner products are not equal!
	\end{example}
	We see that the problem in Example \ref{SameAngleDifferentInnerProducts} is that $\langle \cdot, \cdot \rangle_2$ is a (positive) multiple of that of $\langle \cdot, \cdot \rangle_1$. This phenomenon is called \textbf{\textit{conformality}}, and the two inner products are said to be \textbf{\textit{conformal}}. When two inner products are conformal, the geometric tools (of measuring angles and lengths) are preserved. Particularly, it is easy to see (from a computation similar to that in Example \ref{SameAngleDifferentInnerProducts}) that conformal inner products define the same angles for any pair of vectors. The lengths, however, are scaled by a positive constant. So, if in addition to defining the same angles, if we assume that the two inner products also define the same norm, then we can conclude the equality of the inner products. However, we have remarked earlier that if two inner products define the same norm for \textbf{all} the vectors, then the conclusion follows from the polarization identities (Theorem \ref{PITheorem}), and the hypothesis that they define the same angles becomes redundant! Let us now try to relax the hypothesis. Instead of assuming that the two inner products define the same norms (everywhere), let us see what would happen if there is some $x \in V$ such that $\| x \|_1 = \| x \|_2$. However, in Example \ref{SameAngleDifferentInnerProducts} we again see that the zero vector satisfies this condition but the two inner products are different. So, let us assume that there is some $x \in V \setminus \left\lbrace 0 \right\rbrace$ such that $\| x \|_1 = \| x \|_2$. Can we now say that $\langle \cdot, \cdot \rangle_1 = \langle \cdot, \cdot \rangle_2$? It is clear that if the two inner products are known to be conformal, and there is a non-zero vector with $\| x \|_1 = \| x \|_2$, then the equality of inner products follows. Therefore, we only check what equivalences (in terms of angles) can be drawn with the conformality of the inner products. We divide this analysis into two cases: First we consider the case of real inner product spaces, where things are a bit easy, and then we move to complex inner product spaces.
	\section{The case of real vector spaces}
	\label{RealCase}
\begin{figure}[ht!]
		\centering
		\includegraphics[scale=0.65]{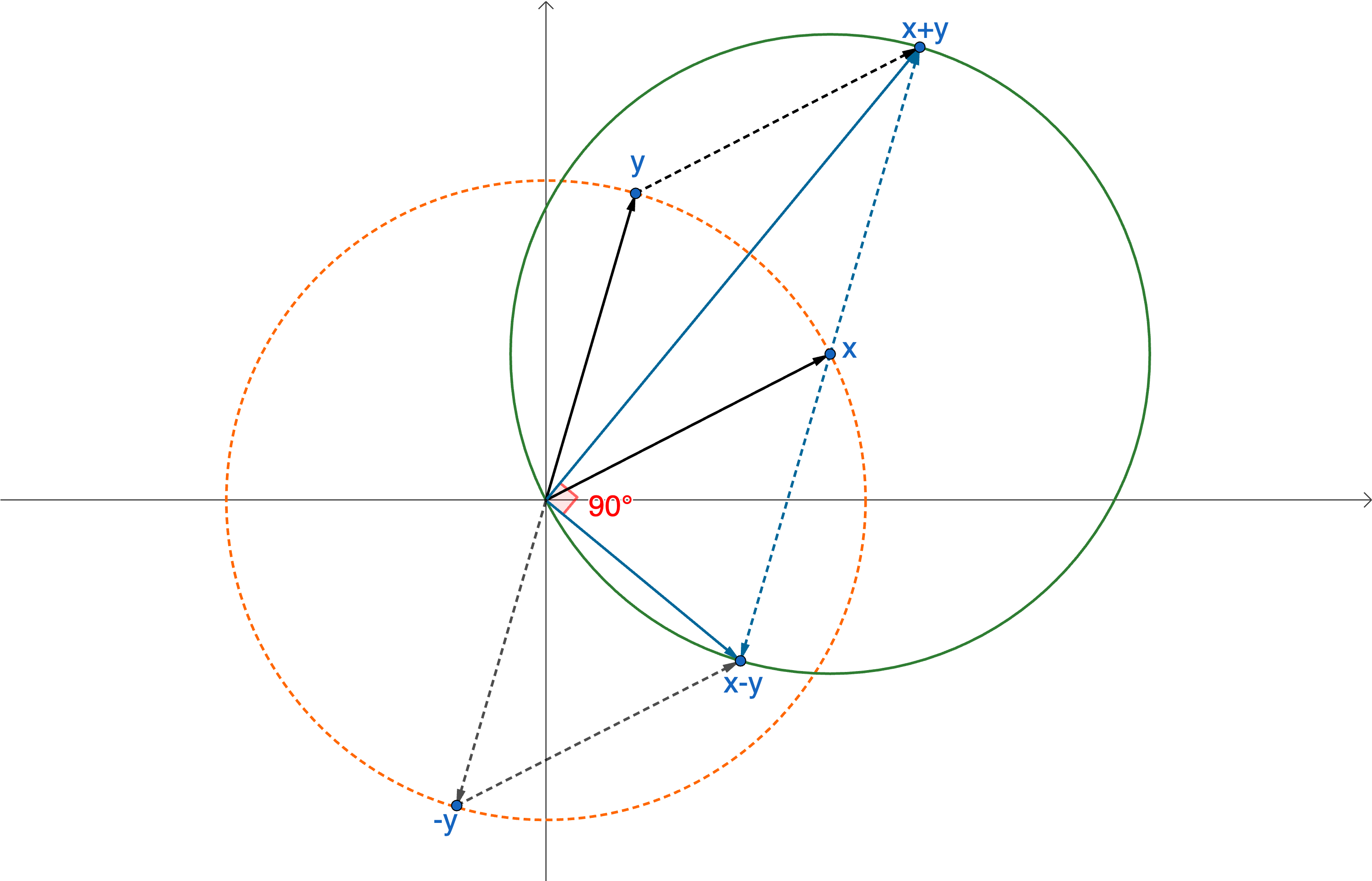}
		\caption{}
		\label{AddSubOrthoFigure}
	\end{figure}

	For this section, we will consider a vector space $V$ over the field of real numbers, $\mathbb{R}$. If $V$ is equipped with only one inner product, we will denote it by $\langle \cdot, \cdot \rangle$, and if it has two inner products, we will denote them by $\langle \cdot, \cdot \rangle_1$ and $\langle \cdot, \cdot \rangle_2$. To answer the question about conformality of the inner products, we begin with the following lemma (see Figure \ref{AddSubOrthoFigure} for a geometric idea).
	\begin{lemma}
		\label{AdditionSubtractionOrthogonal}
		Let $V$ be a real inner product space and let $x, y \in V$ be with $\| x \| = \| y \|$. Then, $\langle x + y, x - y \rangle = 0$.
	\end{lemma}
	\begin{proof}
		From the bilinearity and symmetry of $\langle \cdot, \cdot \rangle$, it easily follows that
		$$\langle x + y, x - y \rangle = \langle x, x \rangle - \langle x, y \rangle + \langle y, x \rangle - \langle y, y \rangle = \| x \|^2 - \| y \|^2 = 0.$$
	\end{proof}
	Now, we are in a position to give a result for real inner product spaces. We will see later in this section that in fact, the result we are about to give can be strengthened. Nonetheless, for starters, the following result\footnote{The statement of this result, given as an exercise problem in the generality of Riemannian metrics in \cite{LeeRM}, has been the motivation for this article.} is an important observation.
	\begin{theorem}
		\label{InnerProductsEqual}
		Let $V$ be a real vector space and $\langle \cdot, \cdot \rangle_1$ and $\langle \cdot, \cdot \rangle_2$ be two inner products on $V$. Assume that the two inner products define the same angles for all pairs of vectors in $V$. Then, there is some $c > 0$ such that $\langle \cdot, \cdot \rangle_2 = c \langle \cdot, \cdot \rangle_1$.
	\end{theorem}
	\begin{proof}
		Let $x, y \in V \setminus \left\lbrace 0 \right\rbrace$. Let $x' = \frac{x}{\| x \|_1}$ and $y' = \frac{y}{\| y \|_1}$. Then, we see that $\| x' \|_1 = \| y' \|_1 = 1$ and hence by Lemma \ref{AdditionSubtractionOrthogonal}, we have $\langle x' + y', x' - y' \rangle_1 = 0$. That is, the angle between $x' + y'$ and $x' - y'$ with respect to $\langle \cdot, \cdot \rangle_1$ is $\frac{\pi}{2}$. Since the two inner products define the same angles, we must also have $\langle x' + y', x' - y' \rangle_2 = 0$. A simple computation then leads to
		\begin{align*}
			\langle x', x' \rangle_2 - \langle y', y' \rangle_2 
			&= 0.
		\end{align*}
		That is, we have for any pair for vectors $x, y \in V \setminus \left\lbrace 0 \right\rbrace$,
		$$\dfrac{\| x \|_2^2}{\| x \|_1^2} = \left\langle \dfrac{x}{\| x \|_1}, \dfrac{x}{\| x \|_1} \right\rangle_2 = \left\langle \dfrac{y}{\| y \|_1}, \dfrac{y}{\| y \|_1} \right\rangle_2 = \dfrac{\| y \|_2^2}{\| y \|_1^2} = c \ (\text{say}).$$
		From Corollary \ref{ConformalCorollary}, we conclude that $\langle \cdot, \cdot \rangle_2 = c \langle \cdot, \cdot \rangle_1$.
	\end{proof}
		We observe that in the proof of Theorem \ref{InnerProductsEqual}, we have used only the fact that the two inner products define the same orthogonality relations ($\frac{\pi}{2}$ angles). Therefore, it is sufficient to state (in the hypothesis) that if two inner products define the same orthogonality relations, then they are conformal. It is also clear, from the following computation that conformal inner products define the same angles.
		$$\cos \theta_2 = \dfrac{\langle x, y \rangle_2}{\| x \|_2 \| y \|_2} = \dfrac{c^2 \langle x, y \rangle_1}{c \| x \|_1 c \| y \|_1} = \dfrac{\langle x, y \rangle_1}{\| x \|_1 \| y \|_1} = \cos \theta_1.$$
		At this point, we may ask what is so special about $\frac{\pi}{2}$ angle? Can we have the result if any $\theta_0$ angle is preserved by the two inner products? Precisely, we raise the following question: \\
		``Let $V$ be a (real) vector space with two inner products $\langle \cdot, \cdot \rangle_1$ and $\langle \cdot, \cdot \rangle_2$, and $\theta_0 \in \left[ 0, \pi \right]$ be fixed. Suppose that the angle between $x, y \in V \setminus \left\lbrace 0 \right\rbrace$ with respect to $\langle \cdot, \cdot \rangle_1$ is $\theta_0$ if and only if the angle between them with respect to $\langle \cdot, \cdot \rangle_2$ is $\theta_0$. Are the two inner products conformal?"

		First we observe that any two inner products define the same $0$ and $\pi$ angles. This follows from the equality in the Cauchy-Schwartz inequality \eqref{CSI}, which is possible if and only if $x$ and $y$ are linearly dependent (see \cite{HoffmannKunze} or \cite{KumaresanLA}). 

\begin{figure}[ht]
			\centering
			\includegraphics[scale=0.65]{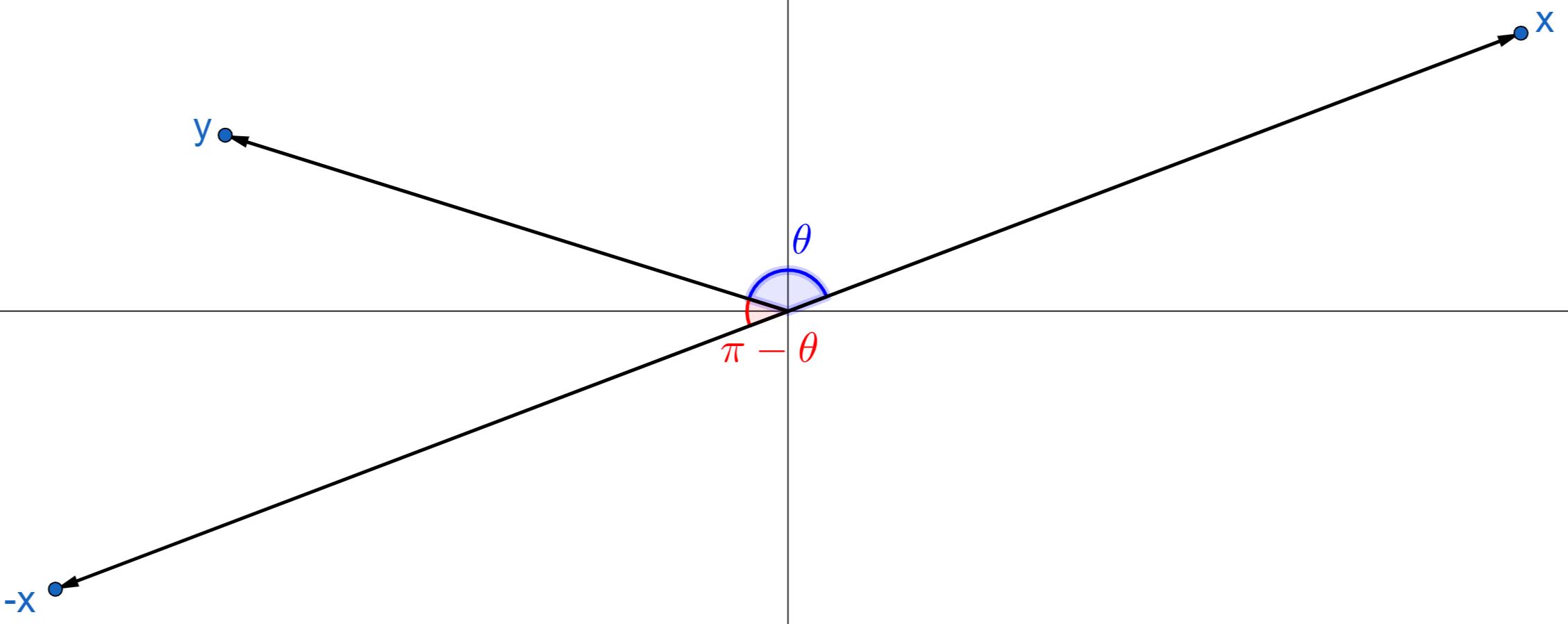}
			\caption{}
			\label{LemmaFigure}
		\end{figure}

		We, therefore, restrict ourselves to $\theta_0 \in \left( 0, \pi \right)$. First we would like to comment that if $\theta_0$ is preserved by two inner products, then so is $\pi - \theta_0$. This is expected because the points $x$ and $-x$ do not depend on the choice of inner product, and they always have angle $\pi$ between them (see Figure \ref{LemmaFigure}).
		\begin{lemma}
			\label{PiMinusThetaNaught}
			Let $V$ be a real vector space with two inner products. Suppose that they define the same $\theta_0$ angle, i.e., for any pair of vectors $x, y \in V$, the angle between $x$ and $y$ with respect to $\langle \cdot, \cdot \rangle_1$ is $\theta_0$ if and only if the angle between them with respect to $\langle \cdot, \cdot \rangle_2$ is $\theta_0$. Then, they define the same $\pi - \theta_0$ angle.
		\end{lemma}
		\begin{proof}
			Let $x, y \in V \setminus \left\lbrace 0 \right\rbrace$ be such that the angle between $x$ and $y$ with respect to $\langle \cdot, \cdot \rangle_1$ is $\pi - \theta_0$. Let $\theta$ be the angle between $-x$ and $y$ with respect to $\langle \cdot, \cdot \rangle_1$. That is,
		$$\cos \theta = \dfrac{\langle -x, y \rangle_1}{\| -x \|_1 \| y \|_1} = - \dfrac{\langle x, y \rangle_1}{\| x \|_1 \| y \|_1} = - \cos \left( \pi - \theta_0 \right) = \cos \theta_0.$$
		Hence, $\theta = \theta_0$. 
		Therefore, the angle between $-x$ and $y$ with respect to $\langle \cdot, \cdot \rangle_2$ is $\theta_0$. Let the angle between $x$ and $y$ with respect to $\langle \cdot, \cdot \rangle_2$ be $\theta$. Then,
		$$\cos \theta = \dfrac{\langle x, y \rangle_2}{\| x \|_2 \| y \|_2} = - \dfrac{\langle -x, y \rangle_2}{\| - x \|_2 \| y \|_2} = - \cos \theta_0.$$
		Hence, we have $\theta = \pi - \theta_0$. This proves the result.
		\end{proof}
		\begin{remark}
			\normalfont
			It is to be noted that Lemma \ref{PiMinusThetaNaught} works both for real and complex vector spaces (except that for complex inner product spaces, we have to use $\text{Re } \langle x, y \rangle$) since we do not use the symmetry (or the bilinearity) of the inner product anywhere in the proof. 
		\end{remark}
		Lemma \ref{PiMinusThetaNaught} significantly reduces our work. Earlier, we would have wanted to prove that whenever two inner products define the same $\theta_0 \in \left( 0, \pi \right)$ angle, they are conformal. Now, we only need to check for $\theta_0 \in \left( 0, \frac{\pi}{2} \right)$, since one of $\theta_0$ or $\pi - \theta_0$ is in $\left( 0, \frac{\pi}{2} \right]$, and Lemma \ref{PiMinusThetaNaught} tells us that these two are preserved simultaneously. Also, we have already discussed the case when $\frac{\pi}{2}$ angle is preserved. We now give the main result of this exposition.
		\begin{theorem}
			\label{SameAngleTheorem}
			Let $V$ be a real vector space equipped with two inner products $\langle \cdot, \cdot \rangle_1$ and $\langle \cdot, \cdot \rangle_2$. Suppose that the two inner products define the same angle $\theta_0 \in \left( 0, \pi \right)$. Then, there is some $c > 0$ such that $\langle \cdot, \cdot \rangle_2 = c \langle \cdot, \cdot \rangle_1$.
		\end{theorem}
		\begin{proof}
			Following the discussion above, we will consider $\theta_0 \in \left( 0, \frac{\pi}{2} \right)$. The advantage of this consideration is that $\sin \theta_0, \cos \theta_0 > 0$. We have seen that to show conformality of inner products it is enough to show that they preserve the $\frac{\pi}{2}$ angle. In this proof, we will employ this technique. \\
			Let $x, y \in V \setminus \left\lbrace 0 \right\rbrace$ be such that $\langle x, y \rangle_1 = 0$. 
			Let the angle between $x$ and $y$ with respect to $\langle \cdot, \cdot \rangle_2$ be $\theta$\footnote{Notice that since $x$ and $y$ are orthogonal in $\langle \cdot, \cdot \rangle_1$, they are linearly independent. Therefore, $\theta \in \left( 0, \pi \right)$.}. Then, we have
			\begin{equation}
				\label{AngleXY2}
				\langle x, y \rangle_2 = \| x \|_2 \| y \|_2 \cos \theta.
			\end{equation}
			Now, we consider the following four vectors as shown in Figure \ref{FourVectors}:
			\begin{align*}
				&z = \cos \theta_0 \dfrac{x}{\| x \|_1} + \sin \theta_0 \dfrac{y}{\| y \|_1}, \
				w = \sin \theta_0 \dfrac{x}{\| x \|_1} + \cos \theta_0 \dfrac{y}{\| y \|_1}, \\
				&\overline{z} = \cos \theta_0 \dfrac{x}{\| x \|_1} - \sin \theta_0 \dfrac{y}{\| y \|_1}, \
				\overline{w} = - \sin \theta_0 \dfrac{x}{\| x \|_1} + \cos \theta_0 \dfrac{y}{\| y \|_1}.
			\end{align*}
			\begin{figure}[ht]
				\centering
				\includegraphics[scale=0.65]{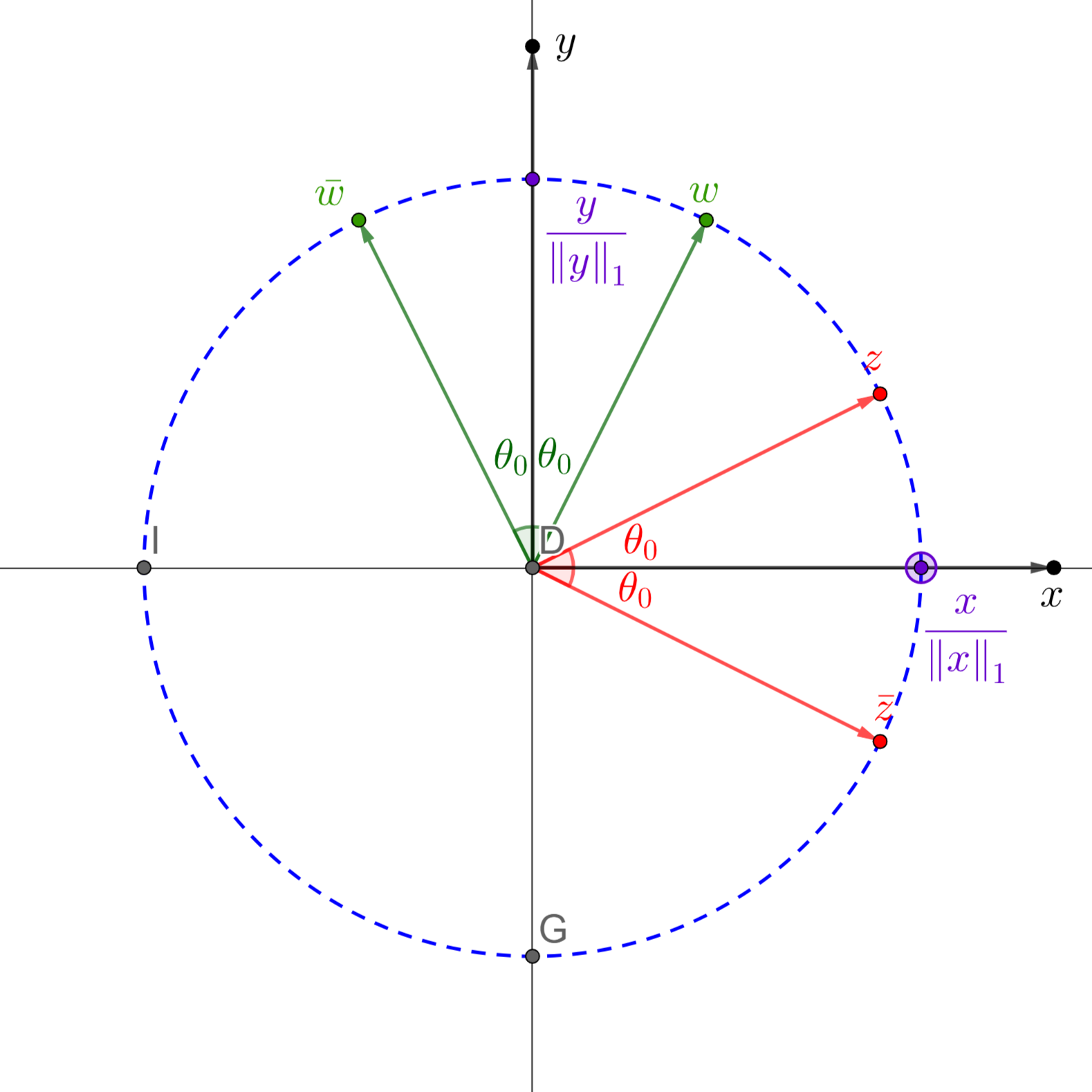}
				\caption{A diagrammatic representation of the vectors $z, w, \overline{z}, \overline{w}$ in the plane spanned by $x$ and $y$.}
				\label{FourVectors}
			\end{figure}
			Since $x$ and $y$ are orthogonal with respect to $\langle \cdot, \cdot \rangle_1$, we immediately see that $\| z \|_1 = \| \bar{z} \|_1 = \| w \|_1 = \| \bar{w} \|_1 = 1$. Moreover, the angle (with respect to $\langle \cdot, \cdot \rangle_1$) between $z$ and $x$, and $\overline{z}$ and $x$ is $\theta_0$, and that between $w$ and $y$, and $\overline{w}$ and $y$ is $\theta_0$. Since the two inner products define the same $\theta_0$ angle, this is true with respect to $\langle \cdot, \cdot \rangle_2$ as well. Thus, we have
		\begin{equation}
			\label{CosThetaNaught}
			\cos \theta_0 = \dfrac{\langle x, z \rangle_2}{\| x \|_2 \| z \|_2} = \dfrac{\langle y, w \rangle_2}{\| y \|_2 \| w \|_2} = \dfrac{\langle x, \overline{z} \rangle_2}{\| x \|_2 \| \overline{z} \|_2} = \dfrac{\langle y, \overline{w} \rangle_2}{\| y \|_2 \| \overline{w} \|_2}.
		\end{equation}
		A simple computation leads us to the following:
		\begin{align}
			\label{NormZSq}
			\| z \|_2^2 &= \cos^2 \theta_0 \dfrac{\| x \|_2^2}{\| x \|_1^2} + 2 \sin \theta_0 \cos \theta_0 \dfrac{\| x \|_2}{\| x \|_1} \dfrac{\| y \|_2}{\| y \|_1} \cos \theta + \sin^2 \theta_0 \dfrac{\| y \|_2^2}{\| y \|_1^2}. \\
			\label{NormWSq}
			\| w \|_2^2 &= \cos^2 \theta_0 \dfrac{\| y \|_2^2}{\| y \|_1^2} + 2 \sin \theta_0 \cos \theta_0 \dfrac{\| x \|_2}{\| x \|_1} \dfrac{\| y \|_2}{\| y \|_1} \cos \theta + \sin^2 \theta_0 \dfrac{\| x \|_2^2}{\| x \|_1^2}. \\
			\label{NormZBarSq}
			\| \overline{z} \|_2^2 &= \cos^2 \theta_0 \dfrac{\| x \|_2^2}{\| x \|_1^2} - 2 \sin \theta_0 \cos \theta_0 \dfrac{\| x \|_2}{\| x \|_1} \dfrac{\| y \|_2}{\| y \|_1} \cos \theta + \sin^2 \theta_0 \dfrac{\| y \|_2^2}{\| y \|_1^2}. \\
			\label{NormWBarSq}
			\| \overline{w} \|_2^2 &= \cos^2 \theta_0 \dfrac{\| y \|_2^2}{\| y \|_1^2} - 2 \sin \theta_0 \cos \theta_0 \dfrac{\| x \|_2}{\| x \|_1} \dfrac{\| y \|_2}{\| y \|_1} \cos \theta + \sin^2 \theta_0 \dfrac{\| x \|_2^2}{\| x \|_1^2}.
		\end{align}
		We also have,
		\begin{align}
			\label{IPXZ2}
			\dfrac{\langle x, z \rangle_2}{\| x \|_2} &= \cos \theta_0 \dfrac{\| x \|_2}{\| x \|_1} + \sin \theta_0 \dfrac{\| y \|_2}{\| y \|_1} \cos \theta. \\
			\label{IPYW2}
			\dfrac{\langle y, w \rangle_2}{\| y \|_2} &= \cos \theta_0 \dfrac{\| y \|_2}{\| y \|_1} + \sin \theta_0 \dfrac{\| x \|_2}{\| x \|_1} \cos \theta. \\
			\label{IPXZBar2}
			\dfrac{\langle x, \overline{z} \rangle_2}{\| x \|_2} &= \cos \theta_0 \dfrac{\| x \|_2}{\| x \|_1} - \sin \theta_0 \dfrac{\| y \|_2}{\| y \|_1} \cos \theta. \\
			\label{IPYWBar2}
			\dfrac{\langle y, \overline{w} \rangle_2}{\| y \|_2} &= \cos \theta_0 \dfrac{\| y \|_2}{\| y \|_1} - \sin \theta_0 \dfrac{\| x \|_2}{\| x \|_1} \cos \theta.
		\end{align}
		First, by squaring Equation \eqref{CosThetaNaught}, 
		using Equations \eqref{NormZSq}, \eqref{NormWSq}, \eqref{IPXZ2}, and \eqref{IPYW2}, and then simplifying, we get
		$$\left( \dfrac{\| y \|_2^2}{\| y \|_1^2} + \dfrac{\| x \|_2^2}{\| x \|_1^2} \right) \left( \dfrac{\| y \|_2^2}{\| y \|_1^2} - \dfrac{\| x \|_2^2}{\| x \|_1^2} \right) = 2 \tan \theta_0 \left( \dfrac{\| x \|_2^2}{\| x \|_1^2} - \dfrac{\| y \|_2^2}{\| y \|_1^2} \right) \dfrac{\| x \|_2}{\| x \|_1} \dfrac{\| y \|_2}{\| y \|_1} \cos \theta.$$
		The final goal of this result is to prove the conformality of the two inner products. If we are successful in proving it, the norms defined by the two inner products must be (positive) multiples of each other. That is to say, we expect (at the end of this proof) that $\frac{\| x \|_2}{\| x \|_1} = \frac{\| y \|_2}{\| y \|_1}$, for all non-zero vectors $x, y \in V$. Let us (for the sake of a contradiction) assume that $\dfrac{\| x \|_2}{\| x \|_1} \neq \dfrac{\| y \|_2}{\| y \|_1}$. 
		Then, we have
		\begin{equation}
			\left( \dfrac{\| y \|_2^2}{\| y \|_1^2} + \dfrac{\| x \|_2^2}{\| x \|_1^2} \right) = - 2 \tan \theta_0 \cos \theta_0 \dfrac{\| x \|_2}{\| x \|_1} \dfrac{\| y \|_2}{\| y \|_1} \cos \theta.
		\end{equation}
		Here, the left hand side is a positive quantity and so are $\sin \theta_0$ and $\cos \theta_0$. This forces $\cos \theta < 0$.

		Again, squaring Equation \eqref{CosThetaNaught}, and proceeding the same way with the Equations \eqref{NormZBarSq}, \eqref{NormWBarSq}, \eqref{IPXZBar2}, and \eqref{IPYWBar2} we get
		\begin{equation}
			\left( \dfrac{\| y \|_2^2}{\| y \|_1^2} + \dfrac{\| x \|_2^2}{\| x \|_1^2} \right) =  2 \tan \theta_0 \cos \theta_0 \dfrac{\| x \|_2}{\| x \|_1} \dfrac{\| y \|_2}{\| y \|_1} \cos \theta,
		\end{equation}
		which forces $\cos \theta > 0$, a contradiction! \\
		Hence, we must have $\dfrac{\| x \|_2}{\| x \|_1} = \dfrac{\| y \|_2}{\| y \|_1}$. Using the first equality of Equation \eqref{CosThetaNaught}, we get
		\begin{align*}
			\cos^2 \theta_0 &= 
			\dfrac{\left( \cos \theta_0 \frac{\| x \|_2}{\| x \|_1} + \sin \theta_0 \frac{\| y \|_2}{\| y \|_1} \cos \theta \right)^2}{\cos^2 \theta_0 \frac{\| x \|_2^2}{\| x \|_1^2} + 2 \sin \theta_0 \cos \theta_0 \frac{\| x \|_2}{\| x \|_1} \frac{\| y \|_2}{\| y \|_1} \cos \theta + \sin^2 \theta_0 \frac{\| y \|_2}{\| y \|_1}} \\\
			&= 1 - \dfrac{\sin^2 \theta_0 \sin^2 \theta}{\cos^2 \theta_0 + 2 \cos \theta_0 \sin \theta_0 \cos \theta + \sin^2 \theta_0}.
		\end{align*}
		Upon simplifying, we get
		\begin{align*}
			\sin^2 \theta \cos^2 \theta_0 &= \left( \cos \theta_0 + \sin \theta_0 \cos \theta \right)^2.
		\end{align*}
		This is same as
		\begin{align*}
			\sin \left( \theta + \theta_0 \right) &= -\cos \theta_0 \text{ or } \sin \left( \theta - \theta_0 \right) = \cos \theta_0.
		\end{align*}
		This gives us that $\theta = \frac{\pi}{2}$ and completes the proof!
		\end{proof}
		We may also ask whether it is possible to ``reduce our work" further. The result we have proved tells us that to check if two inner products are equal, we need to check \textbf{every} pair of vectors that make an angle $\theta_0$ with respect to one of them. Will it be sufficient to only check one pair of vectors with angle $\theta_0$? To put it more precisely, we formulate the following question: \\
		``Let $V$ be a real vector space with two inner products. Suppose that there are $x_0, y_0 \in V \setminus \left\lbrace 0 \right\rbrace$ such that the angle between $x_0$ and $y_0$ in both the inner products is $\theta_0 \in \left( 0, \pi \right)$. Are the two inner products conformal?" 

		The question is negatively answered in the next example.
		\begin{example}
			\label{WeirdGeometryIP}
			\normalfont
			Let us consider $\mathbb{R}^2$. Let $\langle \cdot, \cdot \rangle_1$ be the standard inner product on $\mathbb{R}^2$. Define $\langle \cdot, \cdot \rangle_2 : \mathbb{R}^2 \times \mathbb{R}^2 \rightarrow \mathbb{R}$ as
			\begin{equation}
				\label{WeirdInnerProduct}
				\langle \left( x_1, y_1 \right), \left( x_2, y_2 \right) \rangle_2 = x_1x_2 - \dfrac{1}{2} \left( x_1y_2 + x_2y_1 \right) + \dfrac{1}{2} y_1y_2.
			\end{equation}
			We clearly have,
			$$\| \left( 1, 0 \right) \|_2^2 = 1, \langle \left( 1, 0 \right), \left( 1, 1 \right) \rangle_2 = \dfrac{1}{2}, \text{ and } \| \left( 1, 1 \right) \|_2^2 = \dfrac{1}{2}.$$
			If $\theta$ is the angle between $\left( 1, 0 \right)$ and $\left( 1, 1 \right)$ with respect to $\langle \cdot, \cdot \rangle_2$, then,
			$$\cos \theta = \dfrac{\langle \left( 1, 0 \right), \left( 1, 1 \right) \rangle_2}{\| \left( 1, 0 \right) \|_2 \| \left( 1, 1 \right) \|_2} = \dfrac{1}{\sqrt{2}}.$$
			Thus, $\theta = \frac{\pi}{4}$, which is the same as that with respect to $\langle \cdot, \cdot \rangle_1$. However, it is clear that $\langle \cdot, \cdot \rangle_1 \neq c\langle \cdot, \cdot \rangle_2$, since $c\langle \left( 1, 0 \right), \left( 0, 1 \right) \rangle_2 = - c\frac{1}{2} \neq 0 = \langle \left( 1, 0 \right), \left( 0, 1 \right) \rangle_1$, for any $c > 0$.
			The geometries given by these inner products are shown in Figure \ref{AnglesIP}.
			\begin{figure}[ht]
				\centering
				\begin{minipage}{0.48\textwidth}
					\centering
					\includegraphics[scale=0.55]{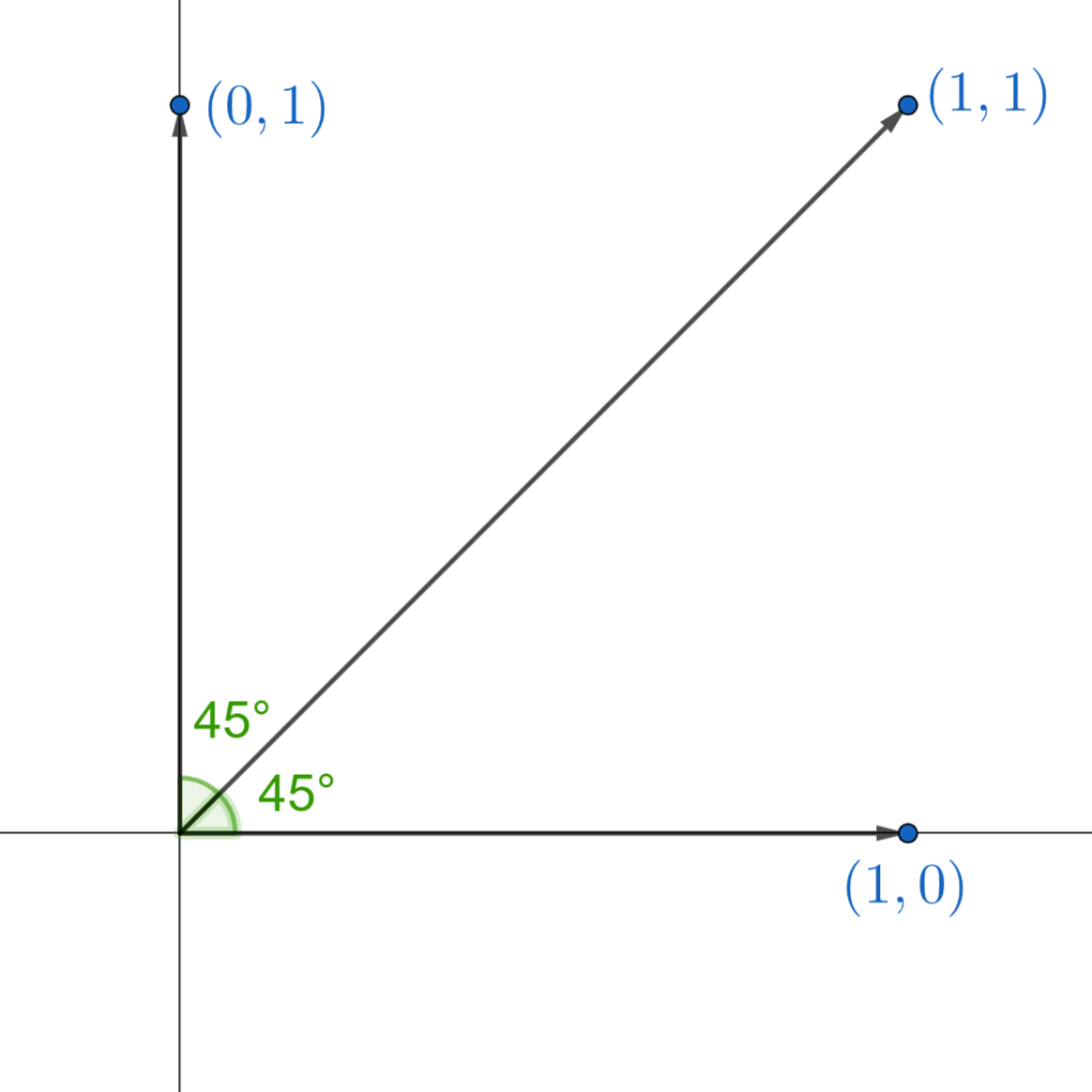}
				\end{minipage}
				\hfill
				\begin{minipage}{0.48\textwidth}
					\centering
					\includegraphics[scale=0.55]{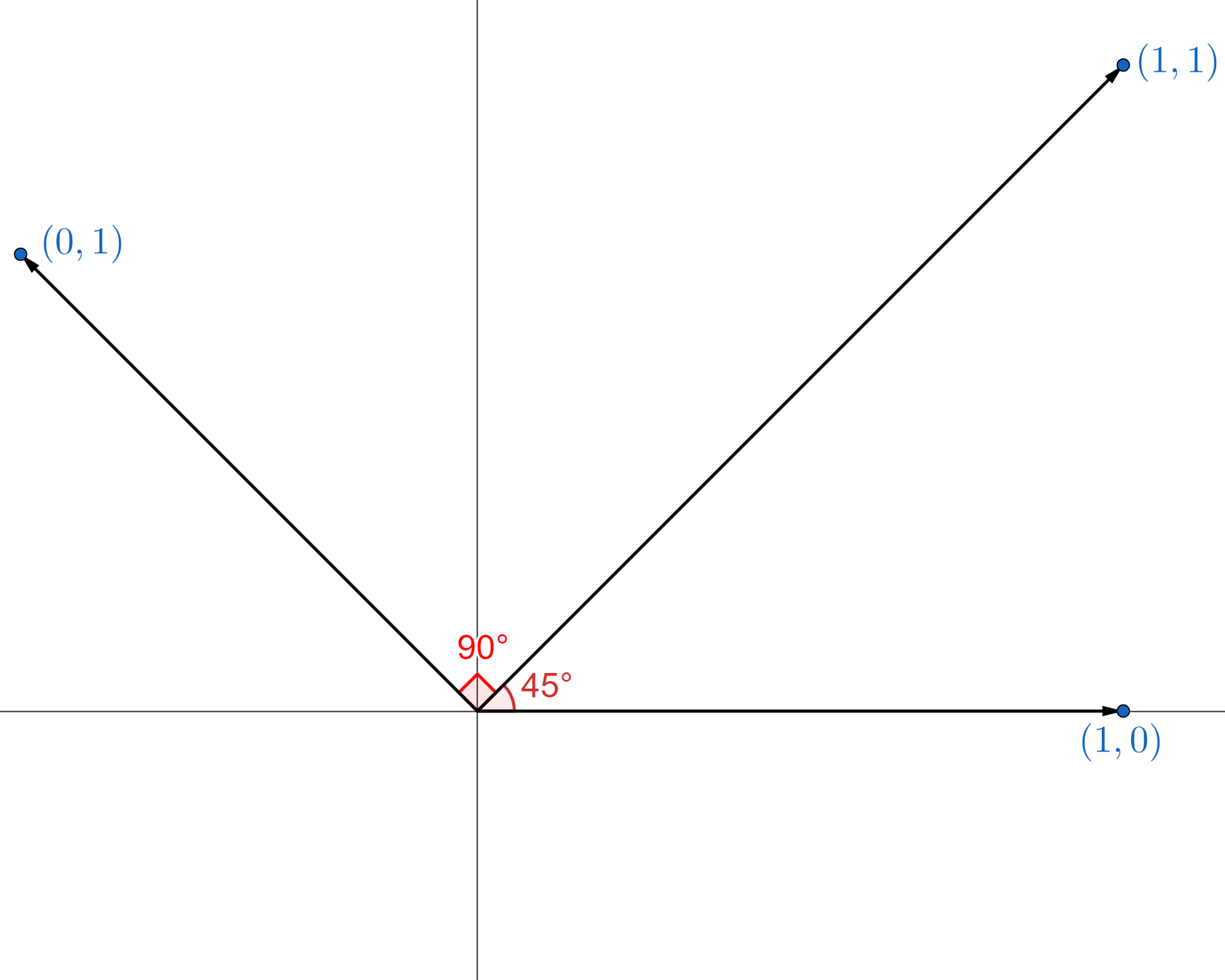}
				\end{minipage}
				\caption{Angles in $\mathbb{R}^2$ in $\langle \cdot, \cdot \rangle_1$ and $\langle \cdot, \cdot \rangle_2$.}
				\label{AnglesIP}
			\end{figure}
		\end{example}
		We summarize the results of this section in the following theorem.
		\begin{theorem}
			\label{SummaryRealCaseTheorem}
			Let $V$ be a real vector space with inner products $\langle \cdot, \cdot \rangle_1$ and $\langle \cdot, \cdot \rangle_2$. Then, the following are equivalent:
			\begin{enumerate}
				\item There is some $c > 0$ such that $\langle \cdot, \cdot \rangle_2 = c \langle \cdot, \cdot \rangle_1$.
				\item There is some $c > 0$ such that $\| \cdot \|_2 = c \| \cdot \|_1$.
				\item The two inner products define the same angles for every pair of vectors of $V$.
				\item The two inner products define the same $\frac{\pi}{2}$ angle (orthogonality relations).
				\item The two inner products define the same $\theta_0 \in \left( 0, \pi \right)$ angle.
			\end{enumerate}
		\end{theorem}
		\begin{corollary}
			\label{SummaryRealCase}
			Let $V$ be a real vector space with inner products $\langle \cdot, \cdot \rangle_1$ and $\langle \cdot, \cdot \rangle_2$. If in addition to any of the conditions in Theorem \ref{SummaryRealCaseTheorem}, there is a non-zero vector $x \in V$ such that $\| x \|_1 = \| x \|_2$, then $\langle \cdot, \cdot \rangle_1 \equiv \langle \cdot, \cdot \rangle_2$.
		\end{corollary}
	\section{The case of complex vector spaces}
	\label{ComplexCase}
	We now ask the question: Can we have results similar to those in Section \ref{RealCase} for complex inner product spaces? We notice that the crux of the proof of Theorem \ref{InnerProductsEqual} was that we wanted to prove $\langle x', x' \rangle_2 - \langle y', y' \rangle_2 = 0$, where $x' = \frac{x}{\| x \|_1}$ and $y' = \frac{y}{\| y \|_1}$, which in turn used Lemma \ref{AdditionSubtractionOrthogonal}. However, Lemma \ref{AdditionSubtractionOrthogonal} will not work in the current state since conjugations will be involved. In complex inner product spaces, due to conjugate symmetry (following the same steps as in Theorem \ref{InnerProductsEqual}), we instead have the following:
	$$\langle x', x' \rangle_2 - \langle y', y' \rangle_2 = \text{Re } \langle x' + y', x' - y' \rangle_2.$$
	Hence, it is now sufficient to show that if $\| x \| = \| y \|$, then $\text{Re } \langle x + y, x - y \rangle = 0$. This easily follows from the following computation,
	\begin{equation}
		\label{AdditionSubtractionPerpendicularEquation}
		\langle x + y, x - y \rangle = \| x \|^2 - \| y \|^2 + 2 \iota \text{Im } \langle y, x \rangle = 2 \iota \text{Im } \langle y, x \rangle.
	\end{equation}
	That is to say, $\text{Re } \langle x + y, x - y \rangle = 0$. Recall that the angle between non-zero vectors $x$ and $y$ in a complex inner product space is given by
	$$\cos \theta = \dfrac{\text{Re } \langle x, y \rangle}{\| x \| \| y \|}.$$
	Thus, the above discussion tells us that for $x, y \in V \setminus \left\lbrace 0 \right\rbrace$ with $\| x \| = \| y \|$, the angle between $x + y$ and $x - y$ is $\frac{\pi}{2}$ (as before!).

	We may ask at this point whether $x + y$ and $x - y$ are orthogonal. For this, we have an easy counterexample on $\mathbb{C}$ (Example \ref{PerpendicularNotOrthogonal}). 

	Using the techniques of Section \ref{RealCase}, we can formulate the analogue of Theorem \ref{SummaryRealCaseTheorem} to complex vector spaces. Particularly, all the equivalences given in Theorem \ref{SummaryRealCaseTheorem} hold verbatim here as well, except that orthogonality relations and the angle $\frac{\pi}{2}$ are no longer equivalent. Let us see if we can have an equivalent criteria for conformality in terms of orthogonality relations as well.

	\begin{theorem}
		\label{OrthogonalityEqualityIP}
		Let $V$ be a complex vector space with two inner products that define the same orthogonality relations. Then, there is some $c > 0$ such that $\langle \cdot, \cdot \rangle_2 = c \langle \cdot, \cdot \rangle_1$.
	\end{theorem}
	\begin{proof}
		Let $x, y \in V \setminus \left\lbrace 0 \right\rbrace$. Let $z = y - \frac{\langle y, x \rangle_1}{\| x \|_1^2} x$. Then, 
		$\langle z, x \rangle_1 = 0$.
		Hence, $\langle z, x \rangle_2 = 0$. Upon substituting the value of $z$ and simplifying, we get
		$\langle y, x \rangle_2 = \langle y, x \rangle_1 \frac{\| x \|_2^2}{\| x \|_1^2}$. 
		Interchanging the roles of $x$ and $y$, we get
		$\langle x, y \rangle_2 = \langle x, y \rangle_1 \frac{\| y \|_2^2}{\| y \|_1^2}$.
		Hence, for all those $x, y \in V \setminus \left\lbrace 0 \right\rbrace$ such that $x$ and $y$ are not orthogonal with respect to $\langle \cdot, \cdot \rangle_1$ and $\langle \cdot, \cdot \rangle_2$, we have
		$\frac{\| x \|_2^2}{\| x \|_1^2} = \frac{\| y \|_2^2}{\| y \|_1^2}$. 
		If $x, y \in V \setminus \left\lbrace 0 \right\rbrace$ are orthogonal in $\langle \cdot, \cdot \rangle_1$. Then, we have $\langle x, y \rangle_1 = \langle x, y \rangle_2 = 0$. Considering $x' = \frac{x}{\| x \|_1}, y' = \frac{y}{\| y \|_1}$, we have from Equation \eqref{AdditionSubtractionPerpendicularEquation}, that $\langle x' + y', x' - y' \rangle_1 = \langle x' + y', x' - y' \rangle_2 = 0$. Upon simplifying, we again get
		$\frac{\| x \|_2^2}{\| x \|_1^2} = \frac{\| y \|_2^2}{\| y \|_1^2} = c \ (\text{say})$.
		The result now follows from Corollary \ref{ConformalCorollary}.
	\end{proof}
	We summarize the results of this section into the following theorem.
	\begin{theorem}
		\label{SummaryComplexCase}
		Let $V$ be a complex vector space with inner products $\langle \cdot, \cdot \rangle_1$ and $\langle \cdot, \cdot \rangle_2$. Then, the following are equivalent:
		\begin{enumerate}
			\item There is a $c > 0$ such that $\langle \cdot, \cdot \rangle_2 = c \langle \cdot, \cdot \rangle_1$.
			\item There is a $c > 0$ such that $\| \cdot \|_2 = c \| \cdot \|_1$.
			\item The two inner products define the same angles for every pair of vectors in $V$.
			\item The two inner products define the same $\frac{\pi}{2}$ angle.
			\item The two inner products define the same orthogonality relations on $V$.
			\item The two inner products define the same $\theta_0 \in \left( 0, \pi \right)$ angle on $V$.
		\end{enumerate}
	\end{theorem}
	\begin{corollary}
		\label{SummaryComplexTheorem}
		Let $V$ be a complex vector space with inner products $\langle \cdot, \cdot \rangle_1$ and $\langle \cdot, \cdot \rangle_2$. If in addition to any of the conditions in Theorem \ref{SummaryComplexCase}, there is a non-zero vector $x \in V$ with $\| x \|_1 = \| x \|_2$, then $\langle \cdot, \cdot \rangle_1 \equiv \langle \cdot, \cdot \rangle_2$.
	\end{corollary}

	\section{Conclusion}
	\label{ConclusionSection}
		In this article, we have discussed some properties of inner products. 
		The main results of the article revolve around knowing when are two inner products defined on a given vector space conformal. What is fascinating about these results is the observation that geometry brings in rigidity. Particularly, we see that if we fix one angle, then all angles are fixed, consequently all projections are scaled and so are the norms. Another fascinating thing about the proofs here is that all these proofs involve only elementary arithmetic and trigonometry.

		We would also like to comment that we can carry forward the results discussed in this article on every tangent space on a Riemannian manifold. Particularly, we have
		\begin{theorem}
			\label{RiemannianManifoldsTheorem}
			Let $M$ be a smooth manifold, and let $g_1, g_2$ be Riemannian metrics on $M$. Then, the following are equivalent:
			\begin{enumerate}
				\item $g_1$ and $g_2$ are (roughly) conformal.
				\item There is a (rough) function $f: M \rightarrow \left( 0, \infty \right)$ such that for each $p \in M$, we have $\| \cdot \|_1^{\left( p \right)} = f \left( p \right) \| \cdot \|_2^{\left( p \right)}$.
				\item For each $p \in M$, $\left( g_1 \right)_p$ and $\left( g_2 \right)_p$ define the same angles on $T_pM$.
				\item For each $p \in M$, $\left( g_1 \right)_p$ and $\left( g_2 \right)_p$ define the same $\frac{\pi}{2}$ angle on $T_pM$.
				\item For each $p \in M$, $\left( g_1 \right)_p$ and $\left( g_2 \right)_p$ define the same $\theta_0 \in \left( 0, \pi \right)$ angle on $T_pM$.
			\end{enumerate}
		\end{theorem}
		It would be interesting to explore if one can get smooth conformality of the Riemannian metrics through smooth vector fields. 
		We end this article with the following question: \\
		
		\textbf{Question:} Given two Riemannian metrics $g_1$ and $g_2$ on a smooth manifold $M$, what are the equivalent criteria for the (smooth) conformality of $g_1$ and $g_2$, preferably in terms of smooth vector fields?
		
\bibliographystyle{amsplain}


\end{document}